\documentclass[11pt]{article}
\usepackage{bbm}
\usepackage{mathrsfs}
\usepackage{amscd}
\usepackage{amsmath,amsfonts,amssymb,amscd}
\usepackage{indentfirst,graphics,epsfig,psfrag}
\input{epsf}
\usepackage{ifpdf}
\usepackage{enumerate}
\usepackage{appendix}
\usepackage{enumerate}
\usepackage{color}
\usepackage{lineno}

\makeatletter \@addtoreset{figure}{section} \makeatother
\makeatletter
\long\def\@makecaption#1#2{%
   \vskip 10\p@
   \setbox\@tempboxa\hbox{{#1}\ \ #2}%
   \ifdim \wd\@tempboxa >\hsize

       {#1}\ \ #2\par
   \else
       \hbox to\hsize{\hfil\box\@tempboxa\hfil}%
   \fi}
\makeatother

\newtheorem{thm}{Theorem}

\newtheorem{lem}{Lemma}

\newtheorem{exe}[thm]{Example}
\newtheorem{obs}{Observation}

\newcommand{\qed}{{\hfill\rule{3pt}{7pt}}}

\def\qed{\hfill \rule{4pt}{7pt}}

\begin{document}
\title{\textbf{Interval minors of complete multipartite graphs}
\footnote{Supported by the National Science Foundation of China
(Nos. 11161037) and the Science Found of Qinghai Province (No.
2014-ZJ-907).}}
\author{
\small Ya-Ping Mao$^1$,  \ Hong-Jian Lai$^2$, \ Zhao Wang$^{1,3}$, \ Zhi-Wei Guo$^1$\\[0.2cm]
\small $^1$Department of Mathematics, Qinghai Normal\\
\small University, Xining, Qinghai 810008, China\\[0.2cm]
\small $^2$Department of Mathematics, West Virginia\\
\small University, Morgantown, WV 26506, USA\\[0.2cm]
\small $^3$School of Mathematical Sciences, Beijing Normal\\
\small University, Beijing 100875, China\\[0.2cm]
\small E-mails: maoyaping@ymail.com; hjlai2015@hotmail.com; \\
\small ~~~~~~~~~wangzhao380@yahoo.com; guozhiweic@yahoo.com}
\date{}
\maketitle
\begin{abstract}
Interval minors of bipartite graphs were introduced by Jacob Fox in
the study of Stanley-Wilf limits. Recently, Mohar, Rafiey,
Tayfeh-Rezaie and Wu investigated the maximum number of edges in
$K_{k,\ell}$-interval minor free bipartite graphs when $k=2$ and
$k=3$. In this paper, we investigate the maximum number of edges in
$K_{k,\ell}$-interval minor free bipartite graphs for general $k$
and $\ell$. We also study the maximum number of edges in
$K_{\ell_1,\ell_2,\cdots,\ell_t}$-interval minor free multipartite
graphs.
\\[2mm]
\\[2mm]
\end{abstract}

\section{Introduction}

All graphs in this paper are undirected, finite and simple. We refer
to \cite{bondy} for undefined graph theoretical notation and terminology.
For a vertex $v$ in a graph $G$, $N_G(v)$ denotes the set of vertices in $G$
adjacent to $v$, called the neighborhood of $v$.
The degree of $v$ in $G$ is $d_G(v) = |N_G(v)|$.
If $X, Y$ are two disjoint vertex subsets of $G$, then $[X,Y]_G$ is the set of
all edges with one end in $X$ and the other end in $Y$.

We follow \cite{Hung74} for the definition of linear orderings of sets.
Throughout this paper, we use
$(A, <_A)$ to denote a linearly ordered set $A$ with a linear ordering $<_A$.
For notational convenience, we often use $A$ to denote $(A, <_A)$ without explicitly
mentioning $<_A$; and when it is clear from the context, we sometimes omit
the subscript $A$ in the linear ordering $<_A$.
Two elements $u$ and $v$
are \emph{consecutive} in the linearly ordered set $A$ if $u<v$ and there is
no vertex $w \in A$ satisfying $u< w<v$.
By \emph{an ordered multipartite graph} $(G;A_1,A_2,\cdots,A_t)$, we mean a $t$-partite graph
$G$ with partite sets $A_1,A_2,\cdots,A_t$ where for each $i$ with $1\leq i \leq t$,
$(A_i, <_{A_i})$ is a linearly ordered set.
All multipartite graphs in this paper are ordered and so, for simplicity, we
usually say multipartite graph $G$ instead of
ordered multipartite graph $(G;A_1,A_2,\cdots,A_t)$.
By identifying two consecutive vertices
$u$ and $v$ to a single vertex $w$ in $G$, we obtain a new ordered bipartite graph
$G'$ such that $N_{G'}(w) = N_G(u) \cup N_G(v)$.

Two ordered bipartite graphs $G$ and $H$ are
\emph{isomorphic} if there is a graph isomorphism $G \rightarrow H$
preserving both parts, possibly exchanging them, and
preserving both linear orders. They are \emph{equivalent} if $H$ can be obtained from
$G$ by reversing the orders in one or both parts of $G$ and possibly exchange
the two parts.

If $G$ and $H$ are ordered bipartite graphs, then $H$ is called an \emph{interval
minor} of $G$ if a graph isomorphic to $H$ can be obtained from $G$ by repeatedly
applying the following operations:
\\
(IM1) deleting an edge;
\\
(IM2) identifying two consecutive vertices.

The operation (IM2) can also be considered as an operation on linearly ordered sets.
Let $(A, <_A)$ be a linearly ordered set with a linear ordering
$a_1 <_A a_2 <_A \cdots, <_A a_n$.
Then for any $i$ with $1 \le i < n$, $a_i, a_{i+1}$ are two consecutive elements in
$A$.  We shall adopt the notational convention of
viewing the operation (IM2) that identifying $a_i$ and $a_{i+1}$ as a mapping
$\phi: A \mapsto A'=A-\{a_{i+1}\}$ by defining $\phi(a_i) = a_i = \phi(a_{i+1})$
and $\phi(a_j) = a_j$ for each $j \neq i, i+1$. Thus $A' = \phi(A)$ has a natural
linear ordering
$a_1 <_{A} a_2 <_{A} \cdots <_{A} a_i <_{A} a_{i+2} <_{A} \cdots <_{A} a_n$
inherited from the linear ordering of $A$. We adopt the convention to view
$A'$ as a linearly ordered subset  of $A$, and to denote this
fact by $A' \subseteq A$.

If $H$ is not an interval minor of $G$, we say that $G$ \emph{avoids}
$H$ as an interval minor or that $G$ is \emph{$H$-interval minor free}.
Let $ex(p,q,H)$ denote the maximum
number of edges in a bipartite graph with parts of sizes $p$ and $q$ which avoids $H$
as an interval minor.

In classical Tur\'{a}n extremal graph theory, one asks about the maximum
number of edges of a graph of order $n$ which has no subgraph isomorphic to a
given graph. Motivated by the problems in computational and combinatorial
geometry, the authors in \cite{BienstockG, Furedi, FurediH}
considered Tur\'{a}n type problems for matrices
which can be seen as ordered bipartite graphs.
In the ordered version
of Tur\'{a}n theory, the question is: what is the maximum number edges of an
ordered bipartite graph with parts of size $p$ and $q$ with no subgraph isomorphic
to a given ordered bipartite graph? For more details on this problem and
its variations, we refer to \cite{AlbertaERWZ, BrassKV, Claessona, Klazar, Mohar, PachT}.
As another variation, interval minors
were recently introduced by Fox in \cite{Fox} in the study of Stanley-Wilf limits. Fox obtained
exponential upper and lower bounds for $ex(n,n,K_{\ell,\ell})$.

Recently, Mohar, Rafiey, Tayfeh-Rezaie and Wu \cite{Mohar}
investigated the maximum number of edges in $K_{k,\ell}$-interval
minor free bipartite graphs when $k=2$ and $\ell=3$. In this paper,
we study the maximum number of edges in $K_{k,\ell}$-interval minor
free bipartite graphs for general $k$ and $\ell$. We also study the
maximum number of edges in
$K_{\ell_1,\ell_2,\cdots,\ell_t}$-interval minor free multipartite
graphs. Our idea is from \cite{Mohar}.

Unless otherwise stated, we in this paper assume that
$0 \le \ell_1\leq \ell_2\leq \cdots \leq \ell_{t}$
are integers.
For notational convenience, we define $m(p,q,k,\ell) = ex(p,q,K_{k,\ell})$, and
$m(n_1,n_2,\cdots,n_{t},\ell_1,\ell_2,\cdots,\ell_{t}) =
ex(n_1,n_2,\cdots,n_{t},K_{\ell_1,\ell_2,\cdots,\ell_{t}})$.
The following observation is immediate.

\begin{obs}\label{obs}
Let $(G;A,B)$ be a bipartite graph with $|A|=p$ and $|B|=q$, and let
$k,\ell$ be two positive integers.

$(1)$ If $K_{k,\ell}$ is an interval minor of $(G;A,B)$, then
$$
\min\{k,\ell\}\leq \min\{p,q\}~and~\max\{k,\ell\}\leq \max\{p,q\}.
$$

$(2)$ If $\min\{p,q\}<\min\{k,\ell\}$ or $\max\{p,q\}<
\max\{k,\ell\}$, then $(G;A,B)$ is $K_{k,\ell}$-interval minor free,
and
$$
m(p,q,k,\ell)=pq.
$$
\end{obs}

The main results are the following theorems, whose proofs are
presented in Sections 2 and 3, respectively.

\begin{thm}\label{th1}
Let $k$ and $\ell$ be two positive integers with $k\leq \ell$,
and let $p$ and $q$ be two positive integers.

$(1)$ If $k\leq p\leq \ell-1$, then
$$
m(p,q,k,\ell)=(\ell-1)(p-k+1)+q(k-1).
$$

$(2)$ If $p=(\ell-k)r+e$, where $k-1\leq e\leq \ell-2$, then
$$
m(p,q,k,\ell)=(\ell-1)(p-k+1)+q(k-1).
$$
\end{thm}

\begin{thm}\label{th2}
Let $n_1,n_2,\cdots,n_t$ be $t$ positive integers, and
$\ell_1,\ell_2,\cdots,\ell_t$ be $t$ positive integers such that
$n_1< n_2< \cdots < n_t$, $\ell_1< \ell_2< \cdots< \ell_t$. If $n_i<\ell_{i+1}$ for $1\leq i\leq t-1$, then
$$
m(n_1,n_2,\cdots,n_{t},\ell_1,\ell_2,\cdots,\ell_{t})= \left(\sum_{i=1}^{t-1}n_i\right)
\left(\sum_{k=i+1}^{\ell}n_k\right)-(\ell_2-1)n_1+(n_2-\ell_2+1)(\ell_1-1).
$$
\end{thm}

\section{Proof of Theorem \ref{th1}}

Throughout this section, we assume that $p,q,k,\ell$ are positive integers.
The purpose of this section is to determine the value of $m(p,q,k,\ell)$
and to complete the proof of Theorem \ref{th1}. Let $(G;A,B)$ be an ordered bipartite graph
where
\begin{eqnarray} \label{order}
& \; & \mbox{ $A$ has a linear ordering
$a_1 <_A a_2 <_A \cdots<_A a_p$ and }
\\ \nonumber
& \; & \mbox{ $B$ has a linear ordering
$b_1<_B b_2 <_B \cdots<_B b_q$.}
\end{eqnarray}

The vertices $a_1$ and $b_1$ are called the \emph{bottom vertices} whereas $a_p$ and
$b_q$ are the \emph{top vertices}.

\begin{lem}\label{lem1}
Suppose that $p \ge k$. Then
$$
m(p,q,k,\ell)\leq (\ell-1)(p-k+1)+q(k-1).
$$
\end{lem}
\begin{pf} Suppose first that $k$ is even. For $k/2\leq i\leq p-k/2$, we define
$X_i=\{b_j \in B\; :$ there exists $i_1<i_2<\cdots<i_{k/2}\leq i< i_{k/2+1}<\cdots<i_{k}$
such that for every $h$ with $1\leq h\leq k$, $a_{i_h}b_j,a_{i}b_j \in E(G) \}$.

If for some $i$, $|X_i| \ge \ell$, then by performing
operations (IM1) and (IM2) to identify vertices in $A$ and deleting the
resulting all but one edge in each resulting parallel class of edges after the vertex identification,
we will obtain a $K_{k,\ell}$-interval minor of $G$, contrary to the assumption.
Hence for every $i$ with $k/2\leq i\leq p-k/2$,  we have $|X_i|\leq \ell-1$.
As there are $k-1$ vertices in $A - \{a_i: k/2\leq i\leq p-k/2\}$,
every $b_j\in B$ appears in at least $d(b_j)-k+1$ sets in $\{X_i: \frac{k}{2}\leq i\leq p-\frac{k}{2}\}$.
Thus
$$
|[B, \cup_{k/2\leq i\leq p-k/2} X_i]_G| \le
\sum_{i=1}^q(d(b_j)-k+1)\leq \sum_{i=\frac{k}{2}}^{p-\frac{k}{2}}|X_i|\leq (p-k+1)(\ell-1).
$$
It follows that $|E(G)| = |[B, \cup_{k/2\leq i\leq p-k/2} X_i]_G|
+ |[B, A-\cup_{k/2\leq i\leq p-k/2} X_i]_G| \leq (p-k+1)(\ell-1)+q(k-1)$.

Suppose that $k$ is odd. For $\frac{k+1}{2}\leq i\leq
p-\frac{k-1}{2}$, we let $X_i=\{b_j \in B\; :$ there exists
$i_1<i_2<\cdots<i_{\frac{k-1}{2}}< i<
i_{\frac{k+1}{2}}<\cdots<i_{k}$ such that for every $h$ with $1\leq
h\leq k$, $a_{i_h}b_j,a_{i}b_j \in E(G) \}$. With a similar argument
as for the case when $k$ is even and by the assumption that $G$ is
$K_{k,\ell}$-interval minor free, we conclude that for each $i$ with
$\frac{k+1}{2}\leq i\leq p-\frac{k-1}{2}$, $|X_i|\leq \ell-1$, and
that every $b_j\in B$ appears in at least $d(b_j)-k+1$ sets in
$\{X_i:\frac{k+1}{2}\leq i\leq p-\frac{k-1}{2}\}$. It follows from
$\displaystyle \sum_{i=1}^q(d(b_j)-k+1)\leq
\sum_{i=\frac{k+1}{2}}^{p-\frac{k-1}{2}}|X_i|\leq (p-k+1)(\ell-1)$.
that $|E(G)|\leq (p-k+1)(\ell-1)+q(k-1)$.\qed
\end{pf}

\begin{exe} \label{pq}
Let $(G;A,B)$ be a bipartite graph with the ordered partite sets $A$ and $B$
as defined in (\ref{order}) and with
\begin{eqnarray} \label{edge-pq}
E(G) & = & \{a_ib_j\,|\,1\leq i\leq p, 1\leq j\leq \ell-1\}
\\ \nonumber
& \; &
\cup \{a_{i_h}b_j\,|\,\ell\leq j\leq q, 1\leq h\leq k-1, i_h\in \{1,2,\cdots,p\}\}.
\end{eqnarray}
(As the edges in $E(G)$ vary as the choice of $i_h$ changes, $(G;A, B)$ defined this way represents
a family of ordered bipartite graphs. We shall use $(G;A,B)$ to denote any one in this family as well. )
\end{exe}

\begin{lem}\label{lem2}
Let $p$ and $q$ be positive integers, and let $k\leq p$ and $\ell\leq q$
be two positive integers with $k\leq p\leq \ell-1$. Then
\\
(i) The graph $(G; A, B)$ defined in Example \ref{pq} is
$K_{k,\ell}$-interval minor free.
\\
(ii) $ m(p,q,k,\ell)\geq (\ell-1)(p-k+1)+q(k-1)$.
\end{lem}
\begin{pf}
Let $(G;A,B)$ be the ordered bipartite graph defined in Example \ref{pq}.
Then direct computation yields  $|E(G)|=(\ell-1)p+(q-\ell+1)(k-1)=(\ell-1)(p-k+1)+q(k-1)$.
Thus (ii) follows from (i).

We argue by contradiction to prove (i) and assume that $G$ has a complete bipartite graph
$(H;A',B')$ as an interval minor with
$k=|A'| \leq |B'|=\ell$. By the definition of an interval
minor, we have either  $A'\subseteq B$ and $B'\subseteq A$,
or $A' \subseteq A$ and $B' \subseteq B$ as linearly
ordered subsets.
If $A'\subseteq B$ and $B'\subseteq A$, then $\ell=|B'|\leq |A|=p$,
contrary to the assumption that $p \leq \ell-1$. Thus
we must have $A'\subseteq
A$ and $B'\subseteq B$.

If $B' \cap \{b_{\ell}, b_{\ell+1}, \cdots, b_q\} \neq \emptyset$,
then there exists a smallest $t$ with $\ell \le t \le q$
such that $b_t \in B'$. By (\ref{edge-pq}), performing (IM2) to identify
consecutive vertices in $B$ will not increase the number of
vertices adjacent to $b_t$, and so
$b_t$ is adjacent to at most $k-1$ vertices in $A'$, contrary to
the fact that $H \cong K_{k,\ell}$. Hence
$B' \subseteq \{b_1, b_2, \cdots, b_{\ell-1}\}$, and so
$|B'| \le \ell-1$, contrary to the assumption that
$|B'| = \ell$. Thus (i) must hold, and so the lemma is justified.
\qed
\end{pf}

Let $(G;A,B)$ and $(G';A',B')$ denote disjoint ordered bipartite
graphs satisfying the following conditions.

$\bullet$ $A$ has ordering $a_1<a_2 <\cdots<a_p$, $A'$ has ordering
$a'_1<a'_2 <\cdots<a'_r$, $B$ has ordering $b_1<b_2 <\cdots< b_q$
and $B'$ has ordering $b'_1<b'_2 <\cdots< b'_t$, where $p,q,r,t\geq
k-1$.

$\bullet$ $a_ib_j\in E(G)$ for each $a_i \ (p-k+2\leq i\leq p)$ and
each $b_j \ (q-k+2\leq j\leq q)$, and $a_i'b_j'\in E(G)$ for each
$a_i' \ (1\leq i\leq k-1)$ and each $b_j' \ (1\leq j\leq k-1)$,
where $a_{p-k+2},a_{p-k+3},\cdots,a_p$ and
$b_{q-k+2},b_{q-k+3},\cdots,b_q$ are the first $k-1$ top vertices of
$A$ and $B$, respectively, and $a'_{1},a'_{2},\cdots,a'_{k-1}$ and
$b'_{1},b'_{2},\cdots,b'_{k-1}$ are the first $k-1$ bottom vertices
of $A'$ and $B'$, respectively.

Denote by $G\oplus G'$ the ordered bipartite graph obtained from
$(G\cup G';A\cup A',B\cup B')$ by identifying $a_{p-k+i}$ with
$a'_{i-1}$ and $b_{q-k+i}$ with $b'_{i-1}$, where $2\leq i\leq k$,
and the linear orders of $A\cup A'$ and $B\cup B'$ are such that the
vertices of $G'$ precede those of $G$. The graph $G\oplus G'$ is
called the \emph{concatenation} of $G$ and $G'$.

In the description of $K_{k,\ell}$-interval minor free graphs below,
we shall use the following simple observation, whose proof is left
to the reader. Let $(G;A,B)$ and $(G';A',B')$ be vertex disjoint
$K_{k,\ell}$-interval minor free bipartite graphs with $k\geq 2$ and
$\ell\geq 2$ such that the $i$-th vertex of the first $k-1$ top
vertices in $A$ and the $i$-th vertex first $k-1$ bottom vertices in
$A'$ are identified to a new vertex, where $1\leq i\leq k-1$, and
the $i$-th vertex of the first $k-1$ top vertices in $B$ and the
$i$-th vertex first $k-1$ bottom vertices in $B'$ are identified to
a new vertex, where $1\leq i\leq k-1$. Then $G\oplus G'$ is also
$K_{k,\ell}$-interval minor free.

\begin{lem}\label{lem3}
Let $p$ and $q$ be positive integers, and let $k$ and $\ell$
be two positive integers with $p=(\ell-k)r+e$, where $k-1\leq e\leq \ell-2$. Then
$$
m(p,q,k,\ell)\geq (\ell-1)(p-k+1)+q(k-1).
$$
\end{lem}
\begin{pf}
We introduce a family of $K_{k,\ell}$-interval minor free bipartite
graphs which would turn out to be extremal. Let $\ell\geq k$ and let $p$ and $q$ be
positive integers and let $r=\lfloor (p-k+1)/(\ell-k)\rfloor$ and
$s=\lfloor (q-k+1)/(\ell-k)\rfloor$. We can write $p=(\ell-k)r+e$ a
nd $q=(\ell-k)s+f$, where $k-1\leq e\leq \ell-2$ and $k-1\leq f\leq \ell-2$.

Suppose now that $r\leq s$. Let $H_0$ be $K_{e,\ell-1}$ and let $H_i$ be a
copy of $K_{\ell-1,\ell-1}$ for $1\leq i\leq r$. The concatenation
$$
H=H_0\oplus H_1\oplus \cdots \oplus H_r.
$$
is $K_{k,\ell}$-interval minor free by the observation preceding
this lemma. It has parts of sizes $p$ and $q'= (\ell-k)(r+1)+(k-1)$.
It also has $r(\ell-k)(\ell+k-2)+e(\ell-1)$ edges. Finally, let
$H^+=K_{k-1,q-q'+(k-1)}$. The graph $H_{p,q}(\ell) = H^+\oplus H$
has parts of sizes $p,q$ and has $(\ell-1)(p-k+1)+q(k-1)$ edges.
Therefore, $m(p,q,k,\ell)\geq (\ell-1)(p-k+1)+q(k-1)$.\qed
\end{pf}

Summing up, Lemmas \ref{lem1} and \ref{lem2} justifies Theorem \ref{th1}(1)
and Lemmas \ref{lem1} and \ref{lem3} justifies Theorem \ref{th1}(2). This completes the proof of Theorem
\ref{th1}.

\section{Proof of Theorem \ref{th2} }

The proof of Theorem \ref{th2} follows immediately
from the following two lemmas.

\begin{lem}\label{lem4}
Let $n_1,n_2,\cdots,n_t$ and $\ell_1,\ell_2,\cdots,\ell_t$ be
positive integers such that $n_1< n_2< \cdots < n_t$, and
$\ell_1< \ell_2< \cdots< \ell_t$. If $n_i<\ell_{i+1}$ for $1\leq i\leq t-1$, then
\begin{eqnarray*}
& \; & m(n_1,n_2,\cdots,n_{t},\ell_1,\ell_2,\cdots,\ell_{t})
\\
& \geq & \left(\sum_{i=1}^{t-1}n_i\right)\left(\sum_{k=i+1}^{\ell}n_k\right)-n_1n_2+(\ell_1-1)n_2+(n_1-\ell_1+1)(\ell_2-1).
\end{eqnarray*}
\end{lem}
\begin{pf}
It suffices to present a complete $t$-partite graph that is
$K_{\ell_1,\ell_2,\cdots,\ell_t}$-interval minor free. Let
$(G;A_1,A_2,\cdots,A_t)$ be a complete $t$-partite graph such that
for each $i$ with  $1\leq i\leq t$, the partite set $A_i$ has a
liner ordering $v_{i,1} < v_{i,2} <\cdots< v_{i,n_i}$; and such that
\begin{eqnarray} \label{edge-2}
E(G)&=&\{v_{1,i}v_{2,j}\,|\,1\leq i\leq n_1, \ 1\leq j\leq \ell_2-1\}
\\ \nonumber
&&\cup \{v_{1,i_h}v_{2,j}\,|\,\ell_2\leq j\leq n_2, \  1\leq h\leq \ell_1-1, \  i_h\in \{1,2,\cdots,n_1\}\}
\\ \nonumber
&&\cup \{v_{i,r}v_{j,s}\,|\,3\leq i,j\leq t, \  1\leq r\leq n_i, \ 1\leq s\leq n_j, \ i\neq j\}
\\ \nonumber
&&\cup \{v_{1,r}v_{j,s}\,|\,3\leq j\leq t, \  1\leq r\leq n_1, \ 1\leq s\leq n_j\}
\\ \nonumber
&&\cup \{v_{2,r}v_{j,s}\,|\,3\leq j\leq t, \  1\leq r\leq n_2, \ 1\leq s\leq n_j\}.
\end{eqnarray}
By (\ref{edge-2}),  $G[A_1\cup A_2]$ is a
complete bipartite graph defined in Example \ref{pq}.
As in  Lemma \ref{lem2}, $(G;A_1,A_2,\cdots,A_t)$
defined this way represents a family of ordered multipartite graphs.
We will also use $(G;A_1,A_2,\cdots,A_t)$ to denote any one in this
family.

We claim that $G$ is $K_{\ell_1,\ell_2,\cdots,\ell_{t}}$-interval
minor free. Assume, to the contrary, that $G$ contains a
$K_{\ell_1,\ell_2,\cdots,\ell_{t}}$-interval minor
$(H;A'_1,A'_2,\cdots,A'_t)$, such that for some
permutation $\tau$ on the set $\{1, 2, \cdots, t\}$,
$A_i'\subseteq A_{\tau(i)}$ as a linearly ordered subset and
$|A_i'| = \ell_i$,  for $1\leq i\leq t$.
For $1\leq i\leq t-1$, since
$n_i<\ell_{i+1}$, it follows that
$A_i'\subseteq A_i$. Since  $H \cong K_{\ell_1,\ell_2,\cdots,\ell_{t}}$
is a $K_{\ell_1,\ell_2,\cdots,\ell_{t}}$-interval minor of $G$.
it follows that $G[A_1\cup A_2]$ contains a complete
bipartite graph $K_{\ell_1,\ell_2}$ as its subgraph,
contrary to  Lemma \ref{lem2}(i). As direct computation yields
\begin{eqnarray*}
|E(G)|&=&\left(\sum_{i=1}^{t-1}n_i\right)\left(\sum_{k=i+1}^{\ell}n_k\right)-n_1n_2+(\ell_2-1)n_1+(n_2-\ell_2+1)(\ell_1-1)\\
&&\left(\sum_{i=1}^{t-1}n_i\right)\left(\sum_{k=i+1}^{\ell}n_k\right)-n_1n_2+(\ell_1-1)n_2+(n_1-\ell_1+1)(\ell_2-1),
\end{eqnarray*}
it follows by definition that
$$
m(n_1,n_2,\cdots,n_{t},\ell_1,\ell_2,\cdots,\ell_{t})\geq \left(\sum_{i=1}^{t-1}n_i\right)\left(\sum_{k=i+1}^{\ell}n_k\right)-n_1n_2+(\ell_1-1)n_2+(n_1-\ell_1+1)(\ell_2-1).
$$\qed
\end{pf}

\begin{lem}\label{lem5}
Let $n_1,n_2,\cdots,n_t$ and $\ell_1,\ell_2,\cdots,\ell_t$ be  positive integers such that $n_1<
n_2< \cdots < n_t$ and $\ell_1< \ell_2< \cdots< \ell_t$. Then
\begin{eqnarray*}
& \; & m(n_1,n_2,\cdots,n_{t},\ell_1,\ell_2,\cdots,\ell_{t})
\\
& \leq  &                                                                              \left(\sum_{i=1}^{t-1}n_i\right)\left(\sum_{k=i+1}^{\ell}n_k\right)-n_1n_2+(\ell_1-1)n_2+(n_1-\ell_1+1)(\ell_2-1).
\end{eqnarray*}
\end{lem}
\begin{pf}
Let $(G;A_1,A_2,\cdots,A_t)$ be a $t$-partite graph such that $G$ is $K_{\ell_1,\ell_2,\cdots,\ell_{t}}$-interval minor free.
Then there exists a bipartite graph $(G;A_i,A_j)$ in $G$ induced by the vertices in $A_i\cup A_j$ such that $(G;A_i,A_j)$ is
$K_{\ell_i,\ell_j}$-interval minor free, where $1\leq i\neq j\leq t$. Without loss of generality, let $(G;A_i,A_j)$ be
$K_{\ell_1,\ell_2}$-interval minor free. By Lemma \ref{lem1}, we have
$$
m(n_1,n_2,\ell_1,\ell_2)\leq (\ell_2-1)(n_1-\ell_1+1)+n_2(\ell_1-1),
$$
and hence
\begin{eqnarray*}
& \; &
m(n_1,n_2,\cdots,n_{t},\ell_1,\ell_2,\cdots,\ell_{t})
\\
& \leq  &                                                                              (\ell_1-1)n_2+(n_1-\ell_1+1)(\ell_2-1)+\left(\sum_{i=1}^{t-1}n_i\right)\left(\sum_{k=i+1}^{\ell}n_k\right)-n_1n_2,
\end{eqnarray*}
as desired.\qed
\end{pf}

\end{document}